\documentstyle[12pt]{article}
\date{January 1991}
\def\cal #1{{\fam2 #1}}

\def\div{\mathop{\rm div}}

\newcommand{\qed}{{\unskip\nobreak\hfill\hbox{ $\Box$}\par}}

\makeatletter
\def\newremark#1{\@ifnextchar[{\@ormrk{#1}}{\@nrmrk{#1}}}

\def\@nrmrk#1#2{%
\@ifnextchar[{\@xnrmrk{#1}{#2}}{\@ynrmrk{#1}{#2}}}

\def\@xnrmrk#1#2[#3]{\expandafter\@ifdefinable\csname #1\endcsname
{\@definecounter{#1}\@addtoreset{#1}{#3}%
\expandafter\xdef\csname the#1\endcsname{\expandafter\noexpand
  \csname the#3\endcsname \@rmrkcountersep \@rmrkcounter{#1}}%
\global\@namedef{#1}{\@rmrk{#1}{#2}}\global\@namedef{end#1}{\@endremark}}}

\def\@ynrmrk#1#2{\expandafter\@ifdefinable\csname #1\endcsname
{\@definecounter{#1}%
\expandafter\xdef\csname the#1\endcsname{\@rmrkcounter{#1}}%
\global\@namedef{#1}{\@rmrk{#1}{#2}}\global\@namedef{end#1}{\@endremark}}}

\def\@ormrk#1[#2]#3{\expandafter\@ifdefinable\csname #1\endcsname
  {\global\@namedef{the#1}{\@nameuse{the#2}}%
\global\@namedef{#1}{\@rmrk{#2}{#3}}%
\global\@namedef{end#1}{\@endremark}}}

\def\@rmrk#1#2{\refstepcounter
    {#1}\@ifnextchar[{\@yrmrk{#1}{#2}}{\@xrmrk{#1}{#2}}}

\def\@xrmrk#1#2{\@beginremark{#2}{\csname the#1\endcsname}\ignorespaces}
\def\@yrmrk#1#2[#3]{\@opargbeginremark{#2}{\csname
       the#1\endcsname}{#3}\ignorespaces}

\def\@rmrkcounter#1{\noexpand\arabic{#1}}
\def\@rmrkcountersep{.}
\ifx\@beginremark\undefined
\def\@beginremark#1#2{\rm \trivlist \item[\hskip \labelsep{\bf #1\ #2}]}
\def\@opargbeginremark#1#2#3{\rm \trivlist
      \item[\hskip \labelsep{\bf #1\ #2\ #3}]}
\fi
\def\@endremark{\endtrivlist}

\def\@xnnrmrk#1{\rm \trivlist \item [\hskip \labelsep {\bf #1.}]\ignorespaces}
\def\@ynnrmrk#1[#2]{\rm \trivlist \item [\hskip \labelsep {\bf #1\ #2.}]%
\ignorespaces}
\@namedef{@endremark*}{\@endremark}
\def\newremark#1{
\@namedef{#1*}{\notnumberedrmrk\csname #1\endcsname}
\@ifnextchar [{\@ormrk {#1}}{\@nrmrk {#1}}
}
\newcommand{\notnumberedrmrk}{\def\@rmrk##1##2{\@ifnextchar [{\@ynnrmrk{##2}}%
{\@xnnrmrk{##2}}}}
\newtheorem{lemma}[subsection]{Lemma}

\newtheorem{theor}[subsection]{Theorem}

\newremark{proof}[subsection]{Proof}
\newremark{remark}[subsection]{Remark}

\def\section{\@startsection {section}{1}{\z@}{-3.5ex plus -1ex minus
 -.2ex}{1.5ex plus .2ex}{\large\bf}}
\def\subsection{\@startsection{subsection}{2}{\z@}{-3.25ex plus -1ex
minus
 -.2ex}{-1em}{\normalsize\bf}}
\makeatother

\begin{document}
\sloppy

\title{ On a Grauert-Riemenschneider
vanishing theorem for Frobenius split varieties in characteristic $p$.
}
\author{V. B. Mehta \and Wilberd van der Kallen}
\maketitle

\section{Introduction}
It is known that the Grauert-Riemenschneider vanishing theorem is not
valid
in characteristic $p$ (\cite{not}). Here we show that it may be
restored in the
presence of a suitable Frobenius splitting. The proof uses
interchanging two
projective limits, one involving iterated Frobenius maps,
cf.\ \cite{Boutot} and
\cite{HartshorneSpeiser}, the other coming from Grothendieck's theorem
on formal
functions. That leads to the following general vanishing theorem which
we then apply in the
situation of the Grauert-Riemenschneider theorem.

\begin{theor}\label{general}
Let $\pi :X\rightarrow Y$ be a proper morphism of schemes of finite type
over a perfect field
of characteristic $p>0$. Let $D$ be a closed
subscheme of $X$ with ideal sheaf
$\cal I$, let  $E$ be a closed
subscheme of $Y$  and let $i\geq 0$ such that\\
1. $D$ contains the geometric points of  $\pi^{-1}E$.\\
2. $R^i\pi_\ast(\cal I)$ vanishes off $E$.\\
3. $X$ is Frobenius split, compatibly with $D$.\\
Then  $R^i\pi_\ast(\cal I)$ vanishes on all of $Y$.
\end{theor}

\begin{theor}[Grauert-Riemenschneider with Frobenius
splitting.]\label{special}
Let $\pi :X\rightarrow Y$ be a proper birational morphism of varieties
in characteristic $p>0$ such that:\\
1. $X$ is non-singular and there is $\sigma\in
H^0(X,K_X^{-1})=H^0(X,c_1(X)) $
such that $\sigma^{p-1}$ splits $X$. (cf.~\cite{MehtaRamanathan}.) \\
2. $D=\div(\sigma)$ contains the exceptional locus of $\pi$ set
theoretically.\\
Then $R^i\pi_\ast K_X=0$ for $i>0$.
\end{theor}

\begin{remark}
It will be clear from the proof that many variations on our
Grauert-Riemenschneider theorem are possible. For instance, one may
replace $D$ by some subdivisor which still contains the exceptional
locus, and thus replace $K_X$
in the conclusion by the new $\cal O_X(-D)$. Similarly, the
birationality
assumption may be weakened, as it is used only to conclude that
condition 2 of
\ref{general} is satisfied.
\end{remark}

\section{Proofs}
\subsection{Proof of \ref{special}.} We assume theorem \ref{general}.
For $E$ we take the image of the exceptional locus. Dualizing  $\sigma$
we get a short exact sequence
$$0\rightarrow K_X \rightarrow \cal O_X \rightarrow \cal O_D
\rightarrow 0,
$$
so $K_X$ may be identified with the ideal sheaf $\cal I$ of $D$. That
$D$ is compatibly split is clear from local computations, cf.~Remark on
page 36 of
\cite{MehtaRamanathan}.\qed

\begin{lemma}
Let $\cdots \rightarrow M_2 \rightarrow M_1 \rightarrow M_0$ be a
projective system of artinian
modules over some ring $R$, with transition maps $f^j_i:M_j\rightarrow
M_i$. If $f^i_0$
is nonzero for all $i$, then the projective limit is nonzero.
\end{lemma}
\begin{proof*}
Put $M^{\rm stab}_i=\bigcap_{j\geq i} f^j_i(M_j)$. Then
 $M^{\rm stab}_i=f^k_i(M_k)$ for $k\gg0$.
So
$$ f^{i+1}_i(M^{\rm stab}_{i+1})=f^{i+1}_if^k_{i+1}(M_k)= f^k_i(M_k)=
M^{\rm stab}_i
$$
 for $k\gg0$. Therefore we have a subsystem $(M^{\rm stab}_i)$ with
 nonzero
surjective maps, whence the result.
\end{proof*}

\subsection{Proof of \ref{general}.}
We argue by contradiction.
We may assume $Y$ is affine, so that  $R^i\pi_\ast(\cal I)$ equals
$H^i(X,\cal I)$. Choose an irreducible component, with generic point $y$ say,
of the support on $Y$ of
$H^i(X,\cal I)$, which we suppose to be nonzero. Observe that $y\in E$.
The Frobenius map $F$ as well as its
splitting act on the exact sequence of sheaves
$$0\rightarrow \cal I\rightarrow \cal O_X\rightarrow \cal
O_D\rightarrow 0.
$$
Therefore
the Frobenius and its iterates act by  split injective endomorphisms,
$p$-linear over
$A=\Gamma (Y,\cal O_Y)$, on $H^i(X,\cal I)$,
and the same remains true after localisation and completion at $y$. Let
$R$ be a regular ring
of the form $L[[X_1,\ldots,X_m]]$ mapping onto $A^{\wedge}_y$, where
$L$ is a field of representatives in
the completed local ring $A^{\wedge}_y$.
In the projective system of artinian modules
$$\cdots R\otimes^{p^r} H^i(X,\cal I)^\wedge_y\rightarrow
 R\otimes^{p^{r-1}} H^i(X,\cal I)^\wedge_y\rightarrow \cdots
$$ all maps towards $R\otimes^{p^{0}} H^i(X,\cal I)^\wedge_y=
H^i(X,\cal I)^\wedge_y $ are nonzero.
Here $R\otimes^{p^r}$ refers to base change along the $r$ times
iterated
Frobenius endomorphism of the regular ring $R$, and the projective
system
is thus the one defining the ``leveling''
$G(H^i(X,\cal I)^\wedge_y)$, in the sense of \cite{HartshorneSpeiser},
of  $H^i(X,\cal I)^\wedge_y$
as an $R$ module. The projective limit is nonzero by the Lemma.
On the other hand, as $R$ is a finite free module over $R$ via $F^r$,
one
may also compute
$$G(H^i(X,\cal I)^\wedge_y)=\lim_{\leftarrow_r} R\otimes^{p^r}
H^i(X,\cal I)^\wedge_y$$ as follows
$$\lim_{\leftarrow_r} R\otimes^{p^r} H^i(X,\cal I)^\wedge_y=
\lim_{\leftarrow_r} R\otimes^{p^r}\lim_{\leftarrow_s} H^i(X_s,\cal
I_s)=$$ $$
\lim_{\leftarrow_r}\lim_{\leftarrow_s} R\otimes^{p^r} H^i(X_s,\cal
I_s)=\lim_{\leftarrow_s}\lim_{\leftarrow_r} R\otimes^{p^r} H^i(X_s,\cal
I_s)
$$ where $X_s$ and $\cal I_s$ are the usual thickenings from
Grothendieck's
theorem on formal functions. But by the Artin-Rees lemma the Frobenius
map acts
nilpotently on  $\cal I_s$, (note that some power of $\cal I$
is contained in the pull back of the ideal sheaf of $E$), so
$\lim_{\leftarrow_r} R\otimes^{p^r} H^i(X_s,\cal I_s)$ vanishes.
But then $G(H^i(X,\cal I)^\wedge_y)$ is both nonzero and zero.\qed

\subsection*{Addresses}
\

\

\noindent
V. B. Mehta, School of Mathematics, \\ Tata Institute of Fundamental
Research, \\ Homi Bhabha Road, Bombay 400005\\
INDIA

\

\noindent
Wilberd van der Kallen, Mathematisch Instituut, Budapestlaan 6,\\
P. O. Box 80.010, 3508 TA Utrecht\\
The Netherlands\\
electronic mail: vdkallen@math.ruu.nl

\end{document}